\documentclass[11pt]{amsart}
\usepackage{amsmath,amssymb}
\newtheorem{theorem}{Theorem}
\newtheorem{proposition}[theorem]{Proposition}
\newtheorem{corollary}[theorem]{Corollary}

\begin{document}

\title[Mean curvature flow]{A sharp bound for the inscribed radius under mean curvature flow}
\author{Simon Brendle}
\address{Department of Mathematics \\ Stanford University \\ Stanford, CA 94305}
\thanks{The author was supported in part by the National Science Foundation under grant DMS-1201924.}
\begin{abstract}
We consider a family of embedded, mean convex hypersurfaces which evolve by the mean curvature flow. It follows from general results of White that the inscribed radius at each point on the surface is at least $\frac{c}{H}$, where $c$ is a constant that depends only on the initial data. Andrews recently gave a new proof of that fact using a direct monotonicity argument. In this paper, we improve this result and show that the inscribed radius is at least $\frac{1}{(1+\delta) \, H}$ at each point where the curvature is large.
\end{abstract}
\maketitle

\section{Introduction}

Let $F: M \times [0,T) \to \mathbb{R}^{n+1}$ be a family of embedded, mean convex hypersurfaces which evolve by the mean curvature flow. To fix notation, let $\mu(x,t)$ denote the reciprocal of the inscribed radius. More precisely, we have 
\[\mu(x,t) = \sup_{y \in M \setminus \{x\}} \frac{2 \, \langle F(x,t)-F(y,t),\nu(x,t) \rangle}{|F(x,t)-F(y,t)|^2}.\] 
Note that $\lambda_1 \leq \hdots \leq \lambda_n \leq \mu$, where the $\lambda_i$ are the principal curvatures. Our main result is a sharp estimate for the function $\mu$ in terms of the mean curvature:

\begin{theorem}
\label{main.theorem}
Let $F: M \times [0,T) \to \mathbb{R}^{n+1}$ be a family of embedded mean convex hypersurfaces which evolve by the mean curvature flow. Then, given any constant $\delta>0$, there exists a constant $C(\delta)$ with the property that $\mu \leq (1+\delta) \, H$ whenever $H \geq C(\delta)$. In other words, the inscribed radius is at least $\frac{1}{(1+\delta) \, H}$ at all points where the curvature is greater than $C(\delta)$.
\end{theorem}

Let now $\rho(x,t)$ the reciprocal of the outer radius. More precisely, we put 
\[\rho(x,t) = \max \bigg \{ \sup_{y \in M \setminus \{x\}} \bigg ( -\frac{2 \, \langle F(x,t)-F(y,t),\nu(x,t) \rangle}{|F(x,t)-F(y,t)|^2} \bigg ),0 \bigg \}.\] 
Note that $-\rho \leq \lambda_1 \leq \hdots \leq \lambda_n$. In particular, a convex surface has $\rho = 0$ and the outer radius is infinite in this case. The following result can be viewed as a refinement of the convexity estimates of Huisken and Sinestrari (cf. \cite{Huisken-Sinestrari1}, \cite{Huisken-Sinestrari2}). 

\begin{theorem}
\label{main.theorem.2}
Let $F: M \times [0,T) \to \mathbb{R}^{n+1}$ be a family of embedded mean convex hypersurfaces which evolve by the mean curvature flow. Then, given any constant $\delta>0$, there exists a constant $C(\delta)$ such that $\rho \leq \delta \, H$ whenever $H \geq C(\delta)$. In other words, the outer radius is at least $\frac{1}{\delta H}$ at all points where the curvature is greater than $C(\delta)$.
\end{theorem}

Let us discuss the background of these results. In a series of papers \cite{White1}, \cite{White2}, \cite{White3}, B.~White achieved several breakthroughs in the understanding of singularities of mean curvature flow for embedded hypersurfaces with positive mean curvature. To explain these results, consider a family of embedded mean convex hypersurfaces $M_t$, $t \in [0,T)$, evolving by mean curvature flow, and consider a sequence of times $t_j \to T$ and any sequence of points $p_j \in M_{t_j}$ such that $H(p_j,t_j) \to \infty$. Brian White's results then imply that, after passing to a subsequence, the rescaled surfaces $H(p_j,t_j) \, (M_{t_j} - p_j)$ converge to a smooth convex hypersurface. This result has several important consequences. Among other things, it implies that  the ratio $\frac{\rho}{H}$ must be very small at points where the curvature is large. Theorem \ref{main.theorem.2} provides an alternative proof of this fact. Moreover, White's results imply that $\frac{\mu}{H} \geq c$ for some constant $c$ that depends only on the initial data. In \cite{Andrews}, B.~Andrews obtained a lower bound for the ratio $\frac{\mu}{H}$ using a direct maximum principle argument (see also \cite{Andrews-Langford-McCoy}). Theorem \ref{main.theorem} improves this to a sharp estimate. Our result is inspired by G.~Huisken's curvature pinching estimates for the mean curvature flow (see \cite{Huisken}, \cite{Huisken-Sinestrari3}) and the corresponding estimates of R.~Hamilton for the Ricci flow (cf. \cite{Hamilton1}, \cite{Hamilton2}, \cite{Hamilton3}).

We next discuss the main steps involved in the proof of Theorem \ref{main.theorem}. In Section \ref{calc}, we derive a differential inequality for the function $\mu$. This step shares some common features with \cite{Andrews}. However, we are able to gain an extra gradient term in the evolution equation. A similar gradient term played a crucial role in our solution of the Lawson Conjecture (cf. \cite{Brendle1}, \cite{Brendle2}, \cite{Brendle3}, \cite{Brendle4}). In Section \ref{auxiliary.inequality}, we estimate the Laplacian $\Delta \mu$ from below. In Section \ref{integral.estimates}, we prove $L^p$-estimates for the function $H^{\sigma-1} \, (\mu - (1+\delta) \, H)$, where $\delta$ is a fixed constant and $\sigma$ is very small. We can then use Stampacchia iteration to show that the function $H^{\sigma-1} \, (\mu - (1+\delta) \, H)$ is uniformly bounded from above for some $\sigma>0$. This last step closely follows the work of Huisken \cite{Huisken} (see also \cite{Huisken-Sinestrari3}, Section 4).

\section{A differential inequality for the function $\mu$}

\label{calc}

We first observe that the set $\{(x,t) \in M \times [0,T): \lambda_n(x,t) < \mu(x,t)\}$ is open. Unfortunately, the function $\mu$ is not smooth in general. However, it turns out that $\mu$ is locally Lipschitz continuous and semi-convex. 

\begin{proposition}
\label{regularity} 
Let us fix a point $(\bar{x},\bar{t})$ such that $\lambda_n(\bar{x},\bar{t}) < \mu(\bar{x},\bar{t})$. Then there exists an open neighborhood $U$ of $\bar{x}$ and a real number $\alpha>0$ such that $\mu$ is Lipschitz continuous and semi-convex on the set $U \times (\bar{t}-\alpha,\bar{t}+\alpha)$.
\end{proposition}

\textbf{Proof.} 
By assumption, we can find a point $\bar{y} \in M \setminus \{\bar{x}\}$ such that 
\[\lambda_n(\bar{x},\bar{t}) < \frac{2 \, \langle F(\bar{x},\bar{t})-F(\bar{y},\bar{t}),\nu(\bar{x},\bar{t}) \rangle}{|F(\bar{x},\bar{t})-F(\bar{y},\bar{t})|^2}.\] 
By continuity, we can find an open set $V$ such that $\bar{x} \in V$, $\bar{y} \notin V$, and 
\[\frac{2 \, \langle F(x,t)-F(y,t),\nu(x,t) \rangle}{|F(x,t)-F(y,t)|^2} < \frac{2 \, \langle F(x,t)-F(\bar{y},t),\nu(x,t) \rangle}{|F(x,t)-F(\bar{y},t)|^2}\] 
for all points $x,y \in V$ and all $t \in (\bar{t}-\alpha,\bar{t}+\alpha)$. We next choose an open neighborhood $U$ of $\bar{x}$ such that $\bar{U} \subset V$. Then 
\[\mu(x,t) = \sup_{y \in M \setminus V} \frac{2 \, \langle F(x,t)-F(y,t),\nu(x,t) \rangle}{|F(x,t)-F(y,t)|^2}\] 
for all points $x \in U$ and all $t \in (\bar{t}-\alpha,\bar{t}+\alpha)$. For each $y \in M \setminus V$, the function $x \mapsto \frac{2 \, \langle F(x,t)-F(y,t),\nu(x,t) \rangle}{|F(x,t)-F(y,t)|^2}$ is a smooth function on $U$. Therefore, $\mu$ is Lipschitz continuous and semi-convex on the set $U$. \\

\begin{corollary} 
\label{regularity.2}
There exists a set $A \subset M \times [0,T)$ of measure zero such that $\mu$ admits a second order Taylor expansion at each point $(\bar{x},\bar{t}) \in \{(x,t) \in M \times [0,T): \lambda_n(x,t) < \mu(x,t)\} \setminus A$. 
\end{corollary}

\textbf{Proof.} 
This follows from a theorem of Alexandrov (cf. \cite{Evans-Gariepy}, Section 6.4). \\

\begin{proposition}
\label{evol}
Let us fix a point $(\bar{x},\bar{t}) \in M \times [0,T)$ such that $\lambda_n(\bar{x},\bar{t}) < \mu(\bar{x},\bar{t})$. Moreover, suppose that $U$ is an open neighborhood of $\bar{x}$ and $\Phi: U \times (\bar{t}-\alpha,\bar{t}] \to \mathbb{R}$ is a smooth function such that $\Phi(\bar{x},\bar{t}) = \mu(\bar{x},\bar{t})$ and $\Phi(x,t) \geq \mu(x,t)$ for all points $(x,t) \in U \times (\bar{t}-\alpha,\bar{t}]$. Then 
\[\frac{\partial \Phi}{\partial t} - \Delta \Phi - |A|^2 \, \Phi + \sum_{i=1}^n \frac{2}{\Phi-\lambda_i} \, (D_i \Phi)^2 \leq 0\] 
at the point $(\bar{x},\bar{t})$.
\end{proposition}

\textbf{Proof.} 
Let us define 
\[Z(x,y,t) = \frac{1}{2} \, \Phi(x,t) \, |F(x,t) - F(y,t)|^2 - \langle F(x,t) - F(y,t),\nu(x,t) \rangle.\] 
Since $\Phi(\bar{x},\bar{t}) = \mu(\bar{x},\bar{t}) > \lambda_n(\bar{x},\bar{t})$, we can find a point $\bar{y}$ such that $\bar{x} \neq \bar{y}$ and $Z(\bar{x},\bar{y},\bar{t}) = 0$. By assumption, we have $Z(x,y,t) \geq 0$ for all points $(x,y,t) \in U \times M \times (\bar{t}-\alpha,\bar{t}]$. This implies 
\begin{align*} 
0 = \frac{\partial Z}{\partial x_i}(\bar{x},\bar{y},\bar{t}) 
&= \frac{1}{2} \, \frac{\partial \Phi}{\partial x_i}(\bar{x},\bar{t}) \, |F(\bar{x},\bar{t}) - F(\bar{y},\bar{t})|^2 \\ 
&+ \Phi(\bar{x},\bar{t}) \, \Big \langle F(\bar{x},\bar{t}) - F(\bar{y},\bar{t}),\frac{\partial F}{\partial x_i}(\bar{x},\bar{t}) \Big \rangle \\ 
&- h_i^k(\bar{x},\bar{t}) \, \Big \langle F(\bar{x},\bar{t}) - F(\bar{y},\bar{t}),\frac{\partial F}{\partial x_k}(\bar{x},\bar{t}) \Big \rangle
\end{align*} 
and 
\begin{align*} 
0 = \frac{\partial Z}{\partial y_i}(\bar{x},\bar{y},\bar{t}) 
&= -\Phi(\bar{x},\bar{t}) \, \Big \langle \frac{\partial F}{\partial y_i}(\bar{y},\bar{t}),F(\bar{x},\bar{t}) - F(\bar{y},\bar{t}) \Big \rangle \\ 
&+ \Big \langle \frac{\partial F}{\partial y_i}(\bar{y},\bar{t}),\nu(\bar{x},\bar{t}) \Big \rangle. 
\end{align*} 
We now analyze the Hessian of $Z$ in spatial direction. To that end, we choose geodesic normal coordinates around $\bar{x}$ such that $h_{ij}(\bar{x},\bar{t})$ is a diagonal matrix. The condition $\frac{\partial Z}{\partial x_i}(\bar{x},\bar{y},\bar{t}) = 0$ implies 
\[\Big \langle F(\bar{x},\bar{t}) - F(\bar{y},\bar{t}),\frac{\partial F}{\partial x_i}(\bar{x},\bar{t}) \Big \rangle = -\frac{1}{2} \, \frac{1}{\Phi(\bar{x},\bar{t})-\lambda_i(\bar{x},\bar{t})} \, \frac{\partial \Phi}{\partial x_i}(\bar{x},\bar{t}) \, |F(\bar{x},\bar{t}) - F(\bar{y},\bar{t})|^2.\] 
Moreover, the condition $\frac{\partial Z}{\partial y_i}(\bar{x},\bar{y},\bar{t}) = 0$ implies that the tangent space at $F(\bar{x},\bar{t})$ is obtained from the tangent space at $F(\bar{y},\bar{t})$ by reflection across the hyperplace orthogonal to $F(\bar{x},\bar{t})-F(\bar{y},\bar{t})$. In particular, the normal vector $\nu(\bar{y},\bar{t})$ is obtained from $\nu(\bar{x},\bar{t})$ by reflection across the hyperplane orthogonal to $F(\bar{x},\bar{t})-F(\bar{y},\bar{t})$. This implies 
\begin{align*} 
\nu(\bar{y},\bar{t}) 
&= \nu(\bar{x},\bar{t}) - 2 \, \frac{\langle F(\bar{x},\bar{t})-F(\bar{y},\bar{t}),\nu(\bar{x},\bar{t}) \rangle}{|F(\bar{x},\bar{t})-F(\bar{y},\bar{t})|} \, \frac{F(\bar{x},\bar{t})-F(\bar{y},\bar{t})}{|F(\bar{x},\bar{t})-F(\bar{y},\bar{t})|} \\ 
&= \nu(\bar{x},\bar{t}) - \Phi(\bar{x},\bar{t}) \, (F(\bar{x},\bar{t})-F(\bar{y},\bar{t})). 
\end{align*}
Let us choose geodesic normal coordinates around $\bar{y}$ such that $\frac{\partial F}{\partial y_i}(\bar{y},\bar{t})$ is obtained from $\frac{\partial F}{\partial x_i}(\bar{x},\bar{t})$ by reflection across the hyperplane orthogonal to $F(\bar{x},\bar{t})-F(\bar{y},\bar{t})$. This gives 
\[\frac{\partial F}{\partial y_i}(\bar{y},\bar{t}) = \frac{\partial F}{\partial x_i}(\bar{x},\bar{t}) - 2 \, \frac{\langle F(\bar{x},\bar{t})-F(\bar{y},\bar{t}),\frac{\partial F}{\partial x_i}(\bar{x},\bar{t}) \rangle}{|F(\bar{x},\bar{t})-F(\bar{y},\bar{t})|} \, \frac{F(\bar{x},\bar{t})-F(\bar{y},\bar{t})}{|F(\bar{x},\bar{t})-F(\bar{y},\bar{t})|}.\] 
In particular, we have 
\[\Big \langle \frac{\partial F}{\partial y_i}(\bar{y},\bar{t}),\frac{\partial F}{\partial x_i}(\bar{x},\bar{t}) \Big \rangle = 1 - 2 \, \frac{\langle F(\bar{x},\bar{t})-F(\bar{y},\bar{t}),\frac{\partial F}{\partial x_i}(\bar{x},\bar{t}) \rangle^2}{|F(\bar{x},\bar{t})-F(\bar{y},\bar{t})|^2}\] 
for $i=1,\hdots,n$. Using the Codazzi equations, we obtain 
\begin{align*} 
\sum_{i=1}^n \frac{\partial^2 Z}{\partial x_i^2}(\bar{x},\bar{y},\bar{t}) 
&= \frac{1}{2} \, \Delta \Phi(\bar{x},\bar{t}) \, |F(\bar{x},\bar{t})-F(\bar{y},\bar{t})|^2 \\ 
&+ 2 \sum_{i=1}^n \frac{\partial \Phi}{\partial x_i}(\bar{x},\bar{t}) \, \Big \langle F(\bar{x},\bar{t})-F(\bar{y},\bar{t}),\frac{\partial F}{\partial x_i}(\bar{x},\bar{t}) \Big \rangle \\ 
&- \sum_{i=1}^n \frac{\partial H}{\partial x_i}(\bar{x},\bar{t}) \, \Big \langle F(\bar{x},\bar{t})-F(\bar{y},\bar{t}),\frac{\partial F}{\partial x_i}(\bar{x},\bar{t}) \Big \rangle \\ 
&- H(\bar{x},\bar{t}) \, \Phi(\bar{x},\bar{t}) \, \langle F(\bar{x},\bar{t})-F(\bar{y},\bar{t}),\nu(\bar{x},\bar{t}) \rangle \\ 
&+ |A(\bar{x},\bar{t})|^2 \, \langle F(\bar{x},\bar{t})-F(\bar{y},\bar{t}),\nu(\bar{x},\bar{t}) \rangle \\ 
&+ n \, \Phi(\bar{x},\bar{t}) - H(\bar{x},\bar{t}), 
\end{align*} 
hence 
\begin{align*} 
&\sum_{i=1}^n \frac{\partial^2 Z}{\partial x_i^2}(\bar{x},\bar{y},\bar{t}) \\ 
&= \frac{1}{2} \, \bigg ( \Delta \Phi(\bar{x},\bar{t}) + |A(\bar{x},\bar{t})|^2 \, \Phi(\bar{x},\bar{t}) \\
&\hspace{20mm} - \sum_{i=1}^n \frac{2}{\Phi(\bar{x},\bar{t})-\lambda_i(\bar{x},\bar{t})} \, \Big ( \frac{\partial \Phi}{\partial x_i}(\bar{x},\bar{t}) \Big )^2 \bigg ) \, |F(\bar{x},\bar{t})-F(\bar{y},\bar{t})|^2 \\ 
&- \sum_{i=1}^n \frac{\partial H}{\partial x_i}(\bar{x},\bar{t}) \, \Big \langle F(\bar{x},\bar{t})-F(\bar{y},\bar{t}),\frac{\partial F}{\partial x_i}(\bar{x},\bar{t}) \Big \rangle \\ 
&- H(\bar{x},\bar{t}) \, \Phi(\bar{x},\bar{t}) \, \langle F(\bar{x},\bar{t})-F(\bar{y},\bar{t}),\nu(\bar{x},\bar{t}) \rangle \\ 
&+ n \, \Phi(\bar{x},\bar{t}) - H(\bar{x},\bar{t}). 
\end{align*} 
Moreover, we have  
\begin{align*} 
&\frac{\partial^2 Z}{\partial x_i \, \partial y_i}(\bar{x},\bar{y},\bar{t}) \\ 
&= -\frac{\partial \Phi}{\partial x_i}(\bar{x},\bar{t}) \, \Big \langle \frac{\partial F}{\partial y_i}(\bar{y},\bar{t}),F(\bar{x},\bar{t})-F(\bar{y},\bar{t}) \Big \rangle \\ 
&- (\Phi(\bar{x},\bar{t})-\lambda_i(\bar{x},\bar{t})) \, \Big \langle \frac{\partial F}{\partial y_i}(\bar{y},\bar{t}),\frac{\partial F}{\partial x_i}(\bar{x},\bar{t}) \Big \rangle \\ 
&= \frac{\partial \Phi}{\partial x_i}(\bar{x},\bar{t}) \, \Big \langle F(\bar{x},\bar{t})-F(\bar{y},\bar{t}),\frac{\partial F}{\partial x_i}(\bar{x},\bar{t}) \Big \rangle \\ 
&- (\Phi(\bar{x},\bar{t})-\lambda_i(\bar{x},\bar{t})) \, \Big ( 1 - 2 \, \frac{\langle F(\bar{x},\bar{t})-F(\bar{y},\bar{t}),\frac{\partial F}{\partial x_i}(\bar{x},\bar{t}) \rangle^2}{|F(\bar{x},\bar{t})-F(\bar{y},\bar{t})|^2} \Big ) \\ 
&= -(\Phi(\bar{x},\bar{t})-\lambda_i(\bar{x},\bar{t})) 
\end{align*} 
and 
\begin{align*} 
\frac{\partial^2 Z}{\partial y_i^2}(\bar{x},\bar{y},\bar{t}) 
&= \Phi(\bar{x},\bar{t}) + h_{ii}(\bar{y},\bar{t}) \, \Phi(\bar{x},\bar{t}) \, \langle \nu(\bar{y},\bar{t}),F(\bar{x},\bar{t})-F(\bar{y},\bar{t}) \rangle \\ 
&- h_{ii}(\bar{y},\bar{t}) \, \langle \nu(\bar{y},\bar{t}),\nu(\bar{x},\bar{t}) \rangle \\ 
&= \Phi(\bar{x},\bar{t}) - h_{ii}(\bar{y},\bar{t}). 
\end{align*} 
In the last step, we have used the identity $\nu(\bar{x},\bar{t}) - \Phi(\bar{x},\bar{t}) \, (F(\bar{x},\bar{t})-F(\bar{y},\bar{t})) = \nu(\bar{y},\bar{t})$. Putting these facts together, we obtain 
\begin{align*} 
&\sum_{i=1}^n \Big ( \frac{\partial^2 Z}{\partial x_i^2}(\bar{x},\bar{y},\bar{t}) + 2 \, \frac{\partial^2 Z}{\partial x_i \, \partial y_i}(\bar{x},\bar{y},\bar{t}) + \frac{\partial^2 Z}{\partial y_i^2}(\bar{x},\bar{y},\bar{t}) \Big ) \\ 
&= \frac{1}{2} \, \bigg ( \Delta \Phi(\bar{x},\bar{t}) + |A(\bar{x},\bar{t})|^2 \, \Phi(\bar{x},\bar{t}) \\
&\hspace{20mm} - \sum_{i=1}^n \frac{2}{\Phi(\bar{x},\bar{t})-\lambda_i(\bar{x},\bar{t})} \, \Big ( \frac{\partial \Phi}{\partial x_i}(\bar{x},\bar{t}) \Big )^2 \bigg ) \, |F(\bar{x},\bar{t})-F(\bar{y},\bar{t})|^2 \\ 
&- \sum_{i=1}^n \frac{\partial H}{\partial x_i}(\bar{x},\bar{t}) \, \Big \langle F(\bar{x},\bar{t})-F(\bar{y},\bar{t}),\frac{\partial F}{\partial x_i}(\bar{x},\bar{t}) \Big \rangle \\ 
&- H(\bar{x},\bar{t}) \, \Phi(\bar{x},\bar{t}) \, \langle F(\bar{x},\bar{t})-F(\bar{y},\bar{t}),\nu(\bar{x},\bar{t}) \rangle \\ 
&+ H(\bar{x},\bar{t}) - H(\bar{y},\bar{t}). 
\end{align*} 
On the other hand, using the identity 
\[\frac{\partial}{\partial t} \nu(\bar{x},\bar{t}) = \sum_{i=1}^n \frac{\partial H}{\partial x_i}(\bar{x},\bar{t}) \, \frac{\partial F}{\partial x_i}(\bar{x},\bar{t}),\] 
we deduce that 
\begin{align*} 
\frac{\partial Z}{\partial t}(\bar{x},\bar{y},\bar{t}) 
&= \frac{1}{2} \, \frac{\partial \Phi}{\partial t}(\bar{x},\bar{t}) \, |F(\bar{x},\bar{t}) - F(\bar{y},\bar{t})|^2 \\ 
&- \Phi(\bar{x},\bar{t}) \, \langle F(\bar{x},\bar{t}) - F(\bar{y},\bar{t}),H(\bar{x},\bar{t}) \, \nu(\bar{x},\bar{t}) - H(\bar{y},\bar{t}) \, \nu(\bar{y},\bar{t}) \rangle \\ 
&+ \langle H(\bar{x},\bar{t}) \, \nu(\bar{x},\bar{t}) - H(\bar{y},\bar{t}) \, \nu(\bar{y},\bar{t}),\nu(\bar{x},\bar{t}) \rangle \\ 
&- \sum_{i=1}^n \frac{\partial H}{\partial x_i}(\bar{x},\bar{t}) \, \Big \langle F(\bar{x},\bar{t}) - F(\bar{y},\bar{t}),\frac{\partial F}{\partial x_i}(\bar{x},\bar{t}) \Big \rangle \\ 
&= \frac{1}{2} \, \frac{\partial \Phi}{\partial t}(\bar{x},\bar{t}) \, |F(\bar{x},\bar{t}) - F(\bar{y},\bar{t})|^2 \\ 
&- \sum_{i=1}^n \frac{\partial H}{\partial x_i}(\bar{x},\bar{t}) \, \Big \langle F(\bar{x},\bar{t}) - F(\bar{y},\bar{t}),\frac{\partial F}{\partial x_i}(\bar{x},\bar{t}) \Big \rangle \\ 
&- H(\bar{x},\bar{t}) \, \Phi(\bar{x},\bar{t}) \, \langle F(\bar{x},\bar{t}) - F(\bar{y},\bar{t}),\nu(\bar{x},\bar{t}) \rangle \\ 
&+ H(\bar{x},\bar{t}) - H(\bar{y},\bar{t}). 
\end{align*} 
Here, we have again used the identity $\nu(\bar{x},\bar{t}) - \Phi(\bar{x},\bar{t}) \, (F(\bar{x},\bar{t})-F(\bar{y},\bar{t})) = \nu(\bar{y},\bar{t})$. Thus, we conclude that 
\begin{align*} 
&\frac{\partial Z}{\partial t}(\bar{x},\bar{y},\bar{t}) - \sum_{i=1}^n \Big ( \frac{\partial^2 Z}{\partial x_i^2}(\bar{x},\bar{y},\bar{t}) + 2 \, \frac{\partial^2 Z}{\partial x_i \, \partial y_i}(\bar{x},\bar{y},\bar{t}) + \frac{\partial^2 Z}{\partial y_i^2}(\bar{x},\bar{y},\bar{t}) \Big ) \\ 
&= \frac{1}{2} \, \bigg ( \frac{\partial \Phi}{\partial t}(\bar{x},\bar{t}) - \Delta \Phi(\bar{x},\bar{t}) - |A(\bar{x},\bar{t})|^2 \, \Phi(\bar{x},\bar{t}) \\ 
&\hspace{20mm} + \sum_{i=1}^n \frac{2}{\Phi(\bar{x},\bar{t})-\lambda_i(\bar{x},\bar{t})} \, \Big ( \frac{\partial \Phi}{\partial x_i}(\bar{x},\bar{t}) \Big )^2 \bigg ) \, |F(\bar{x},\bar{t})-F(\bar{y},\bar{t})|^2. 
\end{align*}
Since the function $Z$ attains a local minimum at the point $(\bar{x},\bar{y},\bar{t})$, we have 
\[\frac{\partial Z}{\partial t}(\bar{x},\bar{y},\bar{t}) - \sum_{i=1}^n \Big ( \frac{\partial^2 Z}{\partial x_i^2}(\bar{x},\bar{y},\bar{t}) + 2 \, \frac{\partial^2 Z}{\partial x_i \, \partial y_i}(\bar{x},\bar{y},\bar{t}) + \frac{\partial^2 Z}{\partial y_i^2}(\bar{x},\bar{y},\bar{t}) \Big ) \leq 0.\] From this, the assertion follows. \\

\begin{corollary} 
\label{distributional.version.of.evol}
The function $\mu$ satisfies 
\[\frac{\partial \mu}{\partial t} \leq \Delta \mu + |A|^2 \, \mu - \sum_{i=1}^n \frac{2}{\mu-\lambda_i} \, (D_i \mu)^2\] 
on the set $\{\lambda_n < \mu\}$, where $\Delta \mu$ is interpreted in the sense of distributions.
\end{corollary}

\textbf{Proof.} 
Let us denote by $\chi_{ij}$ the spatial Hessian of $\mu$ in the sense of Alexandrov. It follows from Proposition \ref{evol} that 
\[\frac{\partial \mu}{\partial t} - \text{\rm tr}(\chi) - |A|^2 \, \mu + \sum_{i=1}^n \frac{2}{\mu-\lambda_i} \, (D_i \mu)^2 \leq 0\]
at each point $(\bar{x},\bar{t}) \in \{(x,t) \in M \times [0,T): \lambda_n(x,t) < \mu(x,t)\} \setminus A$. Since $\mu$ is semi-convex, the Hessian in the sense of Alexandrov is dominated by the distributional second derivative (cf. \cite{Evans-Gariepy}, Section 6.4). In particular, we have $\text{\rm tr}(\chi) \leq \Delta \mu$, where $\Delta \mu$ is interpreted in the sense of distributions. Putting these facts together, the assertion follows.

\section{An auxiliary inequality}

\label{auxiliary.inequality}

In this section, we derive another inequality. This inequality uses in a crucial way the convexity estimates of Huisken and Sinestrari \cite{Huisken-Sinestrari2}. These estimates assert that, given any $\varepsilon > 0$, there exists a constant $K_1(\varepsilon)$ such that $\lambda_1 \geq -\varepsilon \, H - K_1(\varepsilon)$. Moreover, $K_1(\varepsilon)$ depends only on $\varepsilon$ and the initial data.

In the following, we will consider a single hypersurface $M_{\bar{t}}$ for some fixed time $\bar{t}$. We will suppress $\bar{t}$ in the notation, as we will only work with a fixed hypersurface.

\begin{proposition} 
\label{aux}
Fix a point $\bar{x} \in M$ such that $\lambda_n(\bar{x}) < \mu(\bar{x})$. Moreover, suppose that $U$ is an open neighborhood of $\bar{x}$ and $\Phi: U \to \mathbb{R}$ is a smooth function such that $\Phi(\bar{x}) = \mu(\bar{x})$ and $\Phi(x) \geq \mu(x)$ for all $x \in U$. Then we have 
\begin{align*} 
0 &\leq \Delta \Phi + \frac{1}{2} \, |A|^2 \, \Phi - \frac{1}{2} \, H \, \Phi^2 + \frac{1}{2} \, n^3 \, (n\varepsilon \, \Phi + K_1(\varepsilon)) \, \Phi^2 \\ 
&+ \sum_{i=1}^n \frac{1}{\Phi-\lambda_i} \, D_i \Phi \, D_i H + \frac{1}{2} \, \big ( H + n^3 \, (n\varepsilon \, \Phi + K_1(\varepsilon)) \big ) \, \sum_{i=1}^n \frac{1}{(\Phi-\lambda_i)^2} \, (D_i \Phi)^2 
\end{align*} 
at the point $\bar{x}$.
\end{proposition}

\textbf{Proof.} 
As above, we define
\[Z(x,y) = \frac{1}{2} \, \Phi(x) \, |F(x) - F(y)|^2 - \langle F(x) - F(y),\nu(x) \rangle.\] 
Since $\Phi(\bar{x}) = \mu(\bar{x}) > \lambda_n(\bar{x})$, we can find a point $\bar{y}$ such that $\bar{x} \neq \bar{y}$ and $Z(\bar{x},\bar{y}) = 0$. By assumption, we have $Z(x,y) \geq 0$ for all points $(x,y) \in U \times M$. It follows from results in Section \ref{calc} that 
\begin{align*} 
\sum_{i=1}^n \frac{\partial^2 Z}{\partial x_i^2}(\bar{x},\bar{y}) 
&\leq \frac{1}{2} \, \Delta \Phi(\bar{x}) \, |F(\bar{x})-F(\bar{y})|^2 \\ 
&- \sum_{i=1}^n \frac{\partial H}{\partial x_i}(\bar{x}) \, \Big \langle F(\bar{x})-F(\bar{y}),\frac{\partial F}{\partial x_i}(\bar{x}) \Big \rangle \\ 
&- H(\bar{x}) \, \Phi(\bar{x}) \, \langle F(\bar{x})-F(\bar{y}),\nu(\bar{x}) \rangle \\ 
&+ |A(\bar{x})|^2 \, \langle F(\bar{x})-F(\bar{y}),\nu(\bar{x}) \rangle \\ 
&+ n \, \Phi(\bar{x}) - H(\bar{x}). 
\end{align*} 
This inequality can be rewritten as 
\begin{align*} 
&\sum_{i=1}^n \frac{\partial^2 Z}{\partial x_i^2}(\bar{x},\bar{y}) \\ 
&\leq \frac{1}{2} \, \bigg ( \Delta \Phi(\bar{x}) + |A(\bar{x})|^2 \, \Phi(\bar{x}) - H(\bar{x}) \, \Phi(\bar{x})^2 \\ 
&\hspace{20mm} + \sum_{i=1}^n \frac{1}{\Phi(\bar{x})-\lambda_i(\bar{x})} \, \frac{\partial \Phi}{\partial x_i}(\bar{x}) \, \frac{\partial H}{\partial x_i}(\bar{x}) \bigg ) \, |F(\bar{x})-F(\bar{y})|^2 \\ 
&+ n \, \Phi(\bar{x}) - H(\bar{x}). 
\end{align*} 
Moreover, we have 
\[\frac{\partial^2 Z}{\partial x_i \, \partial y_i}(\bar{x},\bar{y}) = -(\Phi(\bar{x})-\lambda_i(\bar{x}))\] 
and 
\[\frac{\partial^2 Z}{\partial y_i^2}(\bar{x},\bar{y}) = \Phi(\bar{x}) - h_{ii}(\bar{y}).\] 
In particular, we have $h_{ii}(\bar{y}) \leq \Phi(\bar{x})$, hence $H(\bar{y}) \leq n \, \Phi(\bar{x})$. Consequently, the convexity estimate of Huisken and Sinestrari \cite{Huisken-Sinestrari2} implies that $h_{ii}(\bar{y}) \geq -\varepsilon \, H(\bar{y}) - K_1(\varepsilon) \geq -n\varepsilon \, \Phi(\bar{x}) - K_1(\varepsilon)$. From this, we deduce that 
\[\frac{\partial^2 Z}{\partial y_i^2}(\bar{x},\bar{y}) \leq \Phi(\bar{x}) + n\varepsilon \, \Phi(\bar{x}) + K_1(\varepsilon).\] 
Thus, we conclude that 
\begin{align*} 
&\sum_{i=1}^n \Big ( \frac{\partial^2 Z}{\partial x_i^2}(\bar{x},\bar{y}) + 2 \, \frac{\Phi(\bar{x})-\lambda_i(\bar{x})}{\Phi(\bar{x})} \, \frac{\partial^2 Z}{\partial x_i \, \partial y_i}(\bar{x},\bar{y}) \\ 
&\hspace{20mm} + \frac{(\Phi(\bar{x})-\lambda_i(\bar{x}))^2}{\Phi(\bar{x})^2} \, \frac{\partial^2 Z}{\partial y_i^2}(\bar{x},\bar{y}) \Big ) \\ 
&\leq \frac{1}{2} \, \bigg ( \Delta \Phi(\bar{x}) + |A(\bar{x})|^2 \, \Phi(\bar{x}) - H(\bar{x}) \, \Phi(\bar{x})^2 \\
&\hspace{20mm} + \sum_{i=1}^n \frac{1}{\Phi(\bar{x})-\lambda_i(\bar{x})} \, \frac{\partial \Phi}{\partial x_i}(\bar{x}) \, \frac{\partial H}{\partial x_i}(\bar{x}) \bigg ) \, |F(\bar{x})-F(\bar{y})|^2 \\ 
&+ n \, \Phi(\bar{x}) - H(\bar{x}) - \sum_{i=1}^n \frac{(\Phi(\bar{x})-\lambda_i(\bar{x}))^2}{\Phi(\bar{x})} \\ 
&+ \sum_{i=1}^n \frac{(\Phi(\bar{x})-\lambda_i(\bar{x}))^2}{\Phi(\bar{x})^2} \, (n\varepsilon \, \Phi(\bar{x}) + K_1(\varepsilon)) \\ 
&\leq \frac{1}{2} \, \bigg ( \Delta \Phi(\bar{x}) + |A(\bar{x})|^2 \, \Phi(\bar{x}) - H(\bar{x}) \, \Phi(\bar{x})^2 \\
&\hspace{20mm} + \sum_{i=1}^n \frac{1}{\Phi(\bar{x})-\lambda_i(\bar{x})} \, \frac{\partial \Phi}{\partial x_i}(\bar{x}) \, \frac{\partial H}{\partial x_i}(\bar{x}) \bigg ) \, |F(\bar{x})-F(\bar{y})|^2 \\ 
&+ H(\bar{x}) - \frac{|A(\bar{x})|^2}{\Phi(\bar{x})} + n^3 \, (n\varepsilon \, \Phi(\bar{x}) + K_1(\varepsilon)). 
\end{align*} 
In the last step, we have used the fact that $0 \leq \Phi(\bar{x}) - \lambda_i(\bar{x}) \leq n \, \Phi(\bar{x})$ for $i=1,\hdots,n$ and $\sum_{i=1}^n \frac{(\Phi(\bar{x})-\lambda_i(\bar{x}))^2}{\Phi(\bar{x})^2} \leq n^3$. We now multiply both sides by $\frac{2}{|F(\bar{x})-F(\bar{y})|^2}$. Using the identity 
\begin{align*} 
\frac{1}{|F(\bar{x})-F(\bar{y})|^2} &= \frac{\langle F(\bar{x})-F(\bar{y}),\nu(\bar{x}) \rangle^2}{|F(\bar{x})-F(\bar{y})|^4} + \sum_{i=1}^n \frac{\langle F(\bar{x})-F(\bar{y}),\frac{\partial F}{\partial x_i}(\bar{x}) \rangle^2}{|F(\bar{x})-F(\bar{y})|^4} \\ 
&= \frac{1}{4} \, \bigg ( \Phi(\bar{x})^2 + \sum_{i=1}^n \frac{1}{(\Phi(\bar{x})-\lambda_i(\bar{x}))^2} \, \Big ( \frac{\partial \Phi}{\partial x_i}(\bar{x}) \Big )^2 \bigg ), 
\end{align*} 
we derive the estimate 
\begin{align*} 
&\frac{2}{|F(\bar{x})-F(\bar{y})|^2} \, \sum_{i=1}^n \Big ( \frac{\partial^2 Z}{\partial x_i^2}(\bar{x},\bar{y}) + 2 \, \frac{\Phi(\bar{x})-\lambda_i(\bar{x})}{\Phi(\bar{x})} \, \frac{\partial^2 Z}{\partial x_i \, \partial y_i}(\bar{x},\bar{y}) \\ 
&\hspace{40mm} + \frac{(\Phi(\bar{x})-\lambda_i(\bar{x}))^2}{\Phi(\bar{x})^2} \, \frac{\partial^2 Z}{\partial y_i^2}(\bar{x},\bar{y}) \Big ) \\ 
&\leq \Delta \Phi(\bar{x}) + |A(\bar{x})|^2 \, \Phi(\bar{x}) - H(\bar{x}) \, \Phi(\bar{x})^2 + \sum_{i=1}^n \frac{1}{\Phi(\bar{x})-\lambda_i(\bar{x})} \, \frac{\partial \Phi}{\partial x_i}(\bar{x}) \, \frac{\partial H}{\partial x_i}(\bar{x}) \\ 
&+ \frac{1}{2} \, \Big ( H(\bar{x}) - \frac{|A(\bar{x})|^2}{\Phi(\bar{x})} + n^3 \, (n\varepsilon \, \Phi(\bar{x}) + K_1(\varepsilon)) \Big ) \\ 
&\hspace{20mm} \cdot \bigg ( \Phi(\bar{x})^2 + \sum_{i=1}^n \frac{1}{(\Phi(\bar{x})-\lambda_i(\bar{x}))^2} \, \Big ( \frac{\partial \Phi}{\partial x_i}(\bar{x}) \Big )^2 \bigg ) \\ 
&= \Delta \Phi(\bar{x}) + \frac{1}{2} \, |A(\bar{x})|^2 \, \Phi(\bar{x}) - \frac{1}{2} \, H(\bar{x}) \, \Phi(\bar{x})^2 + \frac{1}{2} \, n^3 \, (n\varepsilon \, \Phi(\bar{x}) + K_1(\varepsilon)) \, \Phi(\bar{x})^2 \\ 
&+ \sum_{i=1}^n \frac{1}{\Phi(\bar{x})-\lambda_i(\bar{x})} \, \frac{\partial \Phi}{\partial x_i}(\bar{x}) \, \frac{\partial H}{\partial x_i}(\bar{x}) \\ 
&+ \frac{1}{2} \, \Big ( H(\bar{x}) - \frac{|A(\bar{x})|^2}{\Phi(\bar{x})} + n^3 \, (n\varepsilon \, \Phi(\bar{x}) + K_1(\varepsilon)) \Big ) \, \sum_{i=1}^n \frac{1}{(\Phi(\bar{x})-\lambda_i(\bar{x}))^2} \, \Big ( \frac{\partial \Phi}{\partial x_i}(\bar{x}) \Big )^2. 
\end{align*} 
Since the function $Z$ attains a local minimum at the point $(\bar{x},\bar{y})$, we have 
\[\sum_{i=1}^n \Big ( \frac{\partial^2 Z}{\partial x_i^2}(\bar{x},\bar{y}) + 2 \, \frac{\Phi(\bar{x})-\lambda_i(\bar{x})}{\Phi(\bar{x})} \, \frac{\partial^2 Z}{\partial x_i \, \partial y_i}(\bar{x},\bar{y}) + \frac{(\Phi(\bar{x})-\lambda_i(\bar{x}))^2}{\Phi(\bar{x})^2} \, \frac{\partial^2 Z}{\partial y_i^2}(\bar{x},\bar{y}) \Big ) \geq 0.\] 
Putting these facts together, the assertion follows. \\

\begin{corollary} 
\label{distributional.version.of.aux}
Let $\eta$ be a nonnegative test function of class $C^1$ such that $\text{\rm supp} \, \eta \subset \{x \in M: \lambda_n(x) < \mu(x)\}$. Then 
\begin{align*} 
0 &\leq -\int_M \langle \nabla \eta,\nabla \mu \rangle + \frac{1}{2} \int_M \eta \, \big ( |A|^2 \, \mu - H \, \mu^2 + n^3 \, (n\varepsilon \, \mu + K_1(\varepsilon)) \, \mu^2 \big ) \\ 
&+ \int_M \eta \, \sum_{i=1}^n \frac{1}{\mu-\lambda_i} \, D_i \mu \, D_i H \\ 
&+ \frac{1}{2} \int_M \eta \, \big ( H + n^3 \, (n\varepsilon \, \mu + K_1(\varepsilon)) \big ) \, \sum_{i=1}^n \frac{1}{(\mu-\lambda_i)^2} \, (D_i \mu)^2. 
\end{align*} 
\end{corollary}

\textbf{Proof.} 
Let us denote by $\chi_{ij}$ the spatial Hessian of $\mu$ in the sense of Alexandrov. It follows from Proposition \ref{aux} that 
\begin{align*} 
0 &\leq \text{\rm tr}(\chi) + \frac{1}{2} \, |A|^2 \, \mu - \frac{1}{2} \, H \, \mu^2 + \frac{1}{2} \, n^3 \, (n\varepsilon \, \mu + K_1(\varepsilon)) \, \mu^2 \\ 
&+ \sum_{i=1}^n \frac{1}{\mu-\lambda_i} \, D_i \mu \, D_i H + \frac{1}{2} \, \big ( H + n^3 \, (n\varepsilon \, \mu + K_1(\varepsilon)) \big ) \, \sum_{i=1}^n \frac{1}{(\mu-\lambda_i)^2} \, (D_i \mu)^2 
\end{align*} 
at each point in $\{\lambda_n < \mu\}$ where $\mu$ admits a second order Taylor expansion. This implies 
\begin{align*} 
0 &\leq \int_M \eta \, \text{\rm tr}(\chi) + \frac{1}{2} \int_M \eta \, \big ( |A|^2 \, \mu - H \, \mu^2 + n^3 \, (n\varepsilon \, \mu + K_1(\varepsilon)) \, \mu^2 \big ) \\ 
&+ \int_M \eta \, \sum_{i=1}^n \frac{1}{\mu-\lambda_i} \, D_i \mu \, D_i H \\ 
&+ \frac{1}{2} \int_M \eta \, \big ( H + n^3 \, (n\varepsilon \, \mu + K_1(\varepsilon)) \big ) \, \sum_{i=1}^n \frac{1}{(\mu-\lambda_i)^2} \, (D_i \mu)^2. 
\end{align*} 
Moreover, since $\mu$ is semi-convex, we have $\text{\rm tr}(\chi) \leq \Delta \mu$, where $\Delta \mu$ is interpreted in the sense of distributions. This gives 
\[\int_M \eta \, \text{\rm tr}(\chi) \leq -\int_M \langle \nabla \mu,\nabla \eta \rangle.\] 
Putting these facts together, the assertion follows. 

\section{Proof of Theorem \ref{main.theorem}} 

\label{integral.estimates}

Suppose that a number $\delta>0$ is given. By the convexity estimate of Huisken and Sinestrari \cite{Huisken-Sinestrari2}, we can find a constant $K_0$, depending only on $\delta$ and the initial data, such that 
\[(n-1) \, \lambda_1 \geq -\frac{\delta}{2} \, H - K_0 \, \min \{H,1\}.\] 
For each $\sigma \in (0,\frac{1}{2})$, we define 
\[f_\sigma = H^{\sigma-1} \, (\mu - (1+\delta) \, H) - K_0\] 
and 
\[f_{\sigma,+} = \max \{f_\sigma,0\}.\] 
On the set $\{f_\sigma \geq 0\}$, we have 
\[\mu \geq (1+\delta) \, H + K_0 \, H^{1-\sigma} \geq (1+\delta) \, H + K_0 \, \min \{H,1\},\] 
hence 
\[\mu - \lambda_n \geq \sum_{i=1}^{n-1} \lambda_i + \delta \, H + K_0 \, \min \{H,1\} \geq \frac{\delta}{2} \, H.\] 
In particular, we have $\{f_\sigma \geq 0\} \subset \{\lambda_n < \mu\}$. Finally, we can find a constant $\Lambda \geq 1$, depending only on the initial data, such that $\mu \leq \Lambda \, H$ and $|A|^2 \leq \Lambda \, H^2$ at all points in spacetime. \\

\begin{proposition}
\label{Lp.bound}
Given any $\delta>0$, we can find a positive constant $c_0$, depending only on $\delta$ and the initial data, with the following property: if $p \geq \frac{1}{c_0}$ and $\sigma \leq c_0 \, p^{-\frac{1}{2}}$, then we have 
\[\frac{d}{dt} \bigg ( \int_{M_t} f_{\sigma,+}^p \bigg ) \leq C \, \sigma \, p \int_{M_t} f_{\sigma,+}^p + \sigma \, p \, K_0^p \int_{M_t} |A|^2\] 
for almost all $t$. Here, $C$ is a positive constant that depends only on $\delta$ and the initial data, but not on $\sigma$ and $p$.
\end{proposition}

\textbf{Proof.} 
By Corollary \ref{distributional.version.of.evol}, we have 
\[\frac{\partial \mu}{\partial t} - \Delta \mu - |A|^2 \, \mu + \sum_{i=1}^n \frac{2}{\mu-\lambda_i} \, (D_i \mu)^2 \leq 0\] 
on the set $\{\lambda_n < \mu\}$, where $\Delta \mu$ is interpreted in the sense of distributions. A straightforward calculation gives 
\begin{align*} 
&\frac{\partial}{\partial t} f_\sigma - \Delta f_\sigma - 2 \, (1-\sigma) \, \Big \langle \frac{\nabla H}{H},\nabla f_\sigma \Big \rangle \\ 
&+ 2 \, \sum_{i=1}^n \frac{H^{\sigma-1}}{\mu-\lambda_i} \, (D_i \mu)^2 - \sigma \, |A|^2 \, (f_\sigma + K_0) \\ 
&\leq -\sigma \, (1-\sigma) \, H^{\sigma-3} \, (\mu - (1+\delta) \, H) \, |\nabla H|^2 \leq 0 
\end{align*} 
on the set $\{f_\sigma \geq 0\}$, where $\Delta f_\sigma$ is again interpreted in the sense of distributions. This implies 
\begin{align*} 
&\frac{d}{dt} \bigg ( \int_{M_t} f_{\sigma,+}^p \bigg ) \\ 
&\leq -p(p-1) \int_{M_t} f_{\sigma,+}^{p-2} \, |\nabla f_\sigma|^2 + 2 \, (1-\sigma) \, p \int_{M_t} f_{\sigma,+}^{p-1} \, \Big \langle \frac{\nabla H}{H},\nabla f_\sigma \Big \rangle \\ 
&- 2p \int_{M_t} \sum_{i=1}^n \frac{f_{\sigma,+}^{p-1} \, H^{\sigma-1}}{\mu-\lambda_i} \, (D_i \mu)^2 + \sigma \, p \int_{M_t} |A|^2 \, f_{\sigma,+}^{p-1} \, (f_\sigma+K_0) 
\end{align*} 
The last term on the right hand side has an unfavorable sign. To estimate this term, we put $\varepsilon = \frac{\delta}{4n^4\Lambda^2}$. Applying Corollary \ref{distributional.version.of.aux} to the test function $\eta = \frac{f_{\sigma,+}^p}{H}$ gives 
\begin{align*} 
&\frac{1}{2} \int_{M_t} \big ( H \, \mu^2 - |A|^2 \, \mu - n^3 \, (n\varepsilon \, \mu + K_1(\varepsilon)) \, \mu^2 \big ) \, \frac{f_{\sigma,+}^p}{H} \\ 
&\leq -\int_M \Big \langle \nabla \Big ( \frac{f_{\sigma,+}^p}{H} \Big ),\nabla \mu \Big \rangle + \int_M \frac{f_{\sigma,+}^p}{H} \, \sum_{i=1}^n \frac{1}{\mu-\lambda_i} \, D_i \mu \, D_i H \\ 
&+ \frac{1}{2} \int_M \frac{f_{\sigma,+}^p}{H} \, \big ( H + n^3 \, (n\varepsilon \, \mu + K_1(\varepsilon)) \big ) \, \sum_{i=1}^n \frac{1}{(\mu-\lambda_i)^2} \, (D_i \mu)^2 \\ 
&\leq p \int_{M_t} \frac{f_{\sigma,+}^{p-1}}{H} \, |\nabla \mu| \, |\nabla f_\sigma| + C \int_{M_t} \frac{f_{\sigma,+}^p}{H^2} \, |\nabla \mu| \, |\nabla H| + C \int_{M_t} \frac{f_{\sigma,+}^p}{H^2} \, |\nabla \mu|^2. 
\end{align*}
Here, $C$ is a positive constant which depends on $\delta$ and the initial data, but not on $\sigma$ and $p$. On the set $\{f_\sigma \geq 0\}$, we have $\mu \geq (1+\delta) \, H$. Moreover, the convexity estimate of Huisken and Sinestrari implies that $|A|^2 \leq (1+\varepsilon) \, H^2 + K_2(\varepsilon)$. Consequently, we have 
\begin{align*} 
&H \, \mu^2 - |A|^2 \, \mu - n^3 \, (n\varepsilon \, \mu + K_1(\varepsilon)) \, \mu^2 \\ 
&\geq (1+\delta) \, H^2 \, \mu - |A|^2 \, \mu - n^3 \, (n\varepsilon \, \mu + K_1(\varepsilon)) \, \Lambda^2 \, H^2 \\ 
&\geq (\delta-\varepsilon) \, H^2 \, \mu - n^3 \, (n\varepsilon \, \mu + K_1(\varepsilon)) \, \Lambda^2 \, H^2 - K_2(\varepsilon) \, \mu \\ 
&\geq \frac{\delta}{2} \, H^2 \, \mu - C \, H 
\end{align*} 
on the set $\{f_\sigma \geq 0\}$. Therefore, we obtain 
\begin{align*} 
\int_{M_t} |A|^2 \, f_{\sigma,+}^p 
&\leq C \, p \int_{M_t} \frac{f_{\sigma,+}^{p-1}}{H} \, |\nabla \mu| \, |\nabla f_\sigma| + C \int_{M_t} \frac{f_{\sigma,+}^p}{H^2} \, |\nabla \mu| \, |\nabla H| \\ 
&+ C \int_{M_t} \frac{f_{\sigma,+}^p}{H^2} \, |\nabla \mu|^2 + C \int_{M_t} f_{\sigma,+}^p, 
\end{align*} 
where $C$ is a positive constant that depends only on $\delta$ and the initial data. Using the pointwise inequality 
\[f_{\sigma,+}^{p-1} \, (f_\sigma+K_0) \leq 2 \, f_{\sigma,+}^p + K_0^p,\] 
we obtain 
\begin{align*} 
&\int_{M_t} |A|^2 \, f_{\sigma,+}^{p-1} \, (f_\sigma+K_0) \\ 
&\leq C \, p \int_{M_t} \frac{f_{\sigma,+}^{p-1}}{H} \, |\nabla \mu| \, |\nabla f_\sigma| + C \int_{M_t} \frac{f_{\sigma,+}^p}{H^2} \, |\nabla \mu| \, |\nabla H| \\ 
&+ C \int_{M_t} \frac{f_{\sigma,+}^p}{H^2} \, |\nabla \mu|^2 + C \int_{M_t} f_{\sigma,+}^p + K_0^p \int_{M_t} |A|^2, 
\end{align*} 
where $C$ is a positive constant that depends only on $\delta$ and the initial data. Putting these facts together, we conclude that 
\begin{align*} 
&\frac{d}{dt} \bigg ( \int_{M_t} f_{\sigma,+}^p \bigg ) \\ 
&\leq -p(p-1) \int_{M_t} f_{\sigma,+}^{p-2} \, |\nabla f_\sigma|^2 + 2 \, (1-\sigma) \, p \int_{M_t} f_{\sigma,+}^{p-1} \, \Big \langle \frac{\nabla H}{H},\nabla f_\sigma \Big \rangle \\ 
&- 2p \int_{M_t} \sum_{i=1}^n \frac{f_{\sigma,+}^{p-1} \, H^{\sigma-1}}{\mu-\lambda_i} \, (D_i \mu)^2 + C \, \sigma \, p^2 \int_{M_t} \frac{f_{\sigma,+}^{p-1}}{H} \, |\nabla \mu| \, |\nabla f_\sigma| \\ 
&+ C \, \sigma \, p \int_{M_t} \frac{f_{\sigma,+}^p}{H^2} \, |\nabla \mu| \, |\nabla H| + C \, \sigma \, p \int_{M_t} \frac{f_{\sigma,+}^p}{H^2} \, |\nabla \mu|^2 \\ 
&+ C \, \sigma \, p \int_{M_t} f_{\sigma,+}^p + \sigma \, p \, K_0^p \int_{M_t} |A|^2, 
\end{align*} 
where $C$ is a positive constant that depends only on $\delta$ and the initial data. Using the identity 
\[\frac{\nabla H}{H} = \frac{\nabla \mu - H^{1-\sigma} \, \nabla f_\sigma}{(1-\sigma) \, \mu + \sigma \, (1+\delta) \, H},\] 
we obtain 
\[\Big \langle \frac{\nabla H}{H},\nabla f_\sigma \Big \rangle \leq \frac{\langle \nabla \mu,\nabla f_\sigma \rangle}{(1-\sigma) \, \mu + \sigma \, (1+\delta) \, H} \leq C \, H^{-1} \, |\nabla \mu| \, |\nabla f_\sigma|\] 
and 
\[\frac{|\nabla H|}{H} \leq C \, \frac{|\nabla \mu|}{H} + C \, \frac{|\nabla f_\sigma|}{f_{\sigma,+}}.\] 
This implies 
\begin{align*} 
&\frac{d}{dt} \bigg ( \int_{M_t} f_{\sigma,+}^p \bigg ) \\ 
&\leq -p(p-1) \int_{M_t} f_{\sigma,+}^{p-2} \, |\nabla f_\sigma|^2 - 2p \int_{M_t} \sum_{i=1}^n \frac{f_{\sigma,+}^{p-1} \, H^{\sigma-1}}{\mu-\lambda_i} \, (D_i \mu)^2 \\ 
&+ C \, (p + \sigma \, p^2) \int_{M_t} \frac{f_{\sigma,+}^{p-1}}{H} \, |\nabla \mu|  \, |\nabla f_\sigma| + C \, \sigma \, p \int_{M_t} \frac{f_{\sigma,+}^p}{H^2} \, |\nabla \mu|^2 \\ 
&+ C \, \sigma \, p \int_{M_t} f_{\sigma,+}^p + \sigma \, p \, K_0^p \int_{M_t} |A|^2, 
\end{align*} 
where $C$ is a positive constant that depends only on $\delta$ and the initial data.

Therefore, we can find a positive constant $c_0$, depending only on $\delta$ and the initial data, with the following property: if $p \geq \frac{1}{c_0}$ and $\sigma \leq c_0 \, p^{-\frac{1}{2}}$, then we have 
\[\frac{d}{dt} \bigg ( \int_{M_t} f_{\sigma,+}^p \bigg ) \leq C \, \sigma \, p \int_{M_t} f_{\sigma,+}^p + \sigma \, p \, K_0^p \int_{M_t} |A|^2\] 
for almost all $t$. This completes the proof of Proposition \ref{Lp.bound}. \\

Proposition \ref{Lp.bound} implies that 
\[\frac{d}{dt} \bigg ( \int_{M_t} (f_{\sigma,+}^p + \sigma \, p \, K_0^p \, \Lambda) \bigg ) \leq C \, \sigma \, p \int_{M_t} f_{\sigma,+}^p\] 
for almost all $t$. Thus, we conclude that 
\[\int_{M_t} (f_{\sigma,+}^p + \sigma \, p \, K_0^p \, \Lambda) \leq e^{C\sigma p t} \int_{M_0} (f_{\sigma,+}^p + \sigma \, p \, K_0^p \, \Lambda),\] 
hence 
\[\bigg ( \int_{M_t} f_{\sigma,+}^p \bigg )^{\frac{1}{p}} \leq e^{C \sigma t} \, |M_0|^{\frac{1}{p}} \, \Big ( \sup_{M_0} f_{\sigma,+}^p + \sigma \, p \, K_0^p \, \Lambda \Big )^{\frac{1}{p}}.\] 
To complete the proof of Theorem \ref{main.theorem} we need one final ingredient: \\

\begin{proposition}
Let $f_{\sigma,k} = H^{\sigma-1} \, (\mu - (1+\delta) \, H) - k$ and $f_{\sigma,k,+} = \max \{f_{\sigma,k},0\}$. Then 
\begin{align*} 
\frac{d}{dt} \bigg ( \int_{M_t} f_{\sigma,k,+}^p \bigg ) 
&\leq -\frac{1}{2} \, p(p-1) \int_{M_t} f_{\sigma,k,+}^{p-2} \, |\nabla f_{\sigma,k}|^2 \\ 
&+ \sigma \, p \int_{M_t} |A|^2 \, f_{\sigma,k,+}^{p-1} \, (f_{\sigma,k}+k) 
\end{align*}
if $k \geq K_0$ and $p$ is sufficiently large.
\end{proposition} 

\textbf{Proof.} 
Assume that $k \geq K_0$. Using the inequality 
\begin{align*} 
\frac{\partial}{\partial t} f_{\sigma,k} 
&\leq \Delta f_{\sigma,k} + 2 \, (1-\sigma) \, \Big \langle \frac{\nabla H}{H},\nabla f_{\sigma,k} \Big \rangle \\ 
&- 2 \, \sum_{i=1}^n \frac{H^{\sigma-1}}{\mu-\lambda_i} \, (D_i \mu)^2 + \sigma \, |A|^2 \, (f_{\sigma,k} + k), 
\end{align*} 
we obtain 
\begin{align*}
&\frac{d}{dt} \bigg ( \int_{M_t} f_{\sigma,k,+}^p \bigg ) \\ 
&\leq -p(p-1) \int_{M_t} f_{\sigma,k,+}^{p-2} \, |\nabla f_{\sigma,k}|^2 + 2 \, (1-\sigma) \, p \int_{M_t} f_{\sigma,k,+}^{p-1} \, \Big \langle \frac{\nabla H}{H},\nabla f_{\sigma,k} \Big \rangle \\ 
&- 2p \int_{M_t} \sum_{i=1}^n \frac{f_{\sigma,k,+}^{p-1} \, H^{\sigma-1}}{\mu-\lambda_i} \, (D_i \mu)^2 + \sigma \, p \int_{M_t} |A|^2 \, f_{\sigma,k,+}^{p-1} \, (f_{\sigma,k}+k). 
\end{align*} 
As above, we have 
\[\Big \langle \frac{\nabla H}{H},\nabla f_{\sigma,k} \Big \rangle \leq \frac{\langle \nabla \mu,\nabla f_{\sigma,k} \rangle}{(1-\sigma) \, \mu + \sigma \, (1+\delta) \, H} \leq C \, H^{-1} \, |\nabla \mu| \, |\nabla f_{\sigma,k}|,\] 
where $C$ depends only on $\delta$ and the initial data. This gives 
\begin{align*}
&\frac{d}{dt} \bigg ( \int_{M_t} f_{\sigma,k,+}^p \bigg ) \\ 
&\leq -p(p-1) \int_{M_t} f_{\sigma,k,+}^{p-2} \, |\nabla f_{\sigma,k}|^2 + C \, p \int_{M_t} \frac{f_{\sigma,k,+}^{p-1}}{H} \, |\nabla \mu| \, |\nabla f_{\sigma,k}| \\ 
&- 2p \int_{M_t} \sum_{i=1}^n \frac{f_{\sigma,k,+}^{p-1} \, H^{\sigma-1}}{\mu-\lambda_i} \, (D_i \mu)^2 + \sigma \, p \int_{M_t} |A|^2 \, f_{\sigma,k,+}^{p-1} \, (f_{\sigma,k}+k). 
\end{align*} From this, the assertion follows. \\

Following the arguments of Huisken \cite{Huisken}, we can now use the Michael-Simon Sobolev inequality (cf. \cite{Michael-Simon}) and Stampacchia iteration to show that $f_{\sigma,k} \leq 0$, where $\sigma$ and $k$ depend only on $\delta$ and the initial data. From this, we deduce that 
\[\mu \leq (1+2\delta) \, H + B,\] 
where $B$ is a positive constant that depends only on $\delta$ and the initial data. This completes the proof of Theorem \ref{main.theorem}.

\section{Proof of Theorem \ref{main.theorem.2}}

In this Section, we describe the proof of Theorem \ref{main.theorem}. We will need the following auxiliary result:

\begin{proposition}
\label{aux.2}
Let us fix a point $(\bar{x},\bar{t}) \in M \times [0,T)$ such that $\rho(\bar{x},\bar{t}) > 0$ and $\rho(\bar{x},\bar{t}) + \lambda_1(\bar{x},\bar{t}) > 0$. Moreover, suppose that $U \subset M \times [0,T)$ is an open neighborhood of $\bar{x}$ and $\Phi: U \times (\bar{t}-\alpha,\bar{t}] \to \mathbb{R}$ is a smooth function such that $\Phi(\bar{x},\bar{t}) = \rho(\bar{x},\bar{t})$ and $\Phi(x,t) \geq \rho(x,t)$ for all points $(x,t) \in U \times (\bar{t}-\alpha,\bar{t}]$. Then 
\[\frac{\partial \Phi}{\partial t} + \frac{1}{2} \, H \, \Phi^2 - \sum_{i=1}^n \frac{1}{\Phi+\lambda_i} \, D_i \Phi \, D_i H - \frac{1}{2} \, \sum_{i=1}^n \frac{H}{(\Phi+\lambda_i)^2} \, (D_i \Phi)^2 \leq 0\] 
at the point $(\bar{x},\bar{t})$.
\end{proposition}

\textbf{Proof.} 
Let us define 
\[W(x,y,t) = \frac{1}{2} \, \Phi(x,t) \, |F(x,t) - F(y,t)|^2 + \langle F(x,t) - F(y,t),\nu(x,t) \rangle.\] 
Since $\Phi(\bar{x},\bar{t}) = \rho(\bar{x},\bar{t})$, we can find a point $\bar{y}$ such that $\bar{x} \neq \bar{y}$ and $W(\bar{x},\bar{y},\bar{t}) = 0$. By assumption, we have $W(x,y,t) \geq 0$ for all points $(x,y,t) \in U \times M \times (\bar{t}-\alpha,\bar{t}]$. As in Section \ref{calc}, we compute 
\begin{align*} 
0 = \frac{\partial W}{\partial x_i}(\bar{x},\bar{y},\bar{t}) 
&= \frac{1}{2} \, \frac{\partial \Phi}{\partial x_i}(\bar{x},\bar{t}) \, |F(\bar{x},\bar{t}) - F(\bar{y},\bar{t})|^2 \\ 
&+ \Phi(\bar{x},\bar{t}) \, \Big \langle F(\bar{x},\bar{t}) - F(\bar{y},\bar{t}),\frac{\partial F}{\partial x_i}(\bar{x},\bar{t}) \Big \rangle \\ 
&+ h_i^k(\bar{x},\bar{t}) \, \Big \langle F(\bar{x},\bar{t}) - F(\bar{y},\bar{t}),\frac{\partial F}{\partial x_k}(\bar{x},\bar{t}) \Big \rangle.
\end{align*} 
As usual, we choose geodesic normal coordinates around $\bar{x}$ such that $h_{ij}(\bar{x},\bar{t})$ is a diagonal matrix. The condition $\frac{\partial Z}{\partial x_i}(\bar{x},\bar{y},\bar{t}) = 0$ implies 
\[\Big \langle F(\bar{x},\bar{t}) - F(\bar{y},\bar{t}),\frac{\partial F}{\partial x_i}(\bar{x},\bar{t}) \Big \rangle = -\frac{1}{2} \, \frac{1}{\Phi(\bar{x},\bar{t})+\lambda_i(\bar{x},\bar{t})} \, \frac{\partial \Phi}{\partial x_i}(\bar{x},\bar{t}) \, |F(\bar{x},\bar{t}) - F(\bar{y},\bar{t})|^2.\] 
Moreover, we have 
\begin{align*} 
&\frac{\partial W}{\partial t}(\bar{x},\bar{y},\bar{t}) \\ 
&= \frac{1}{2} \, \frac{\partial \Phi}{\partial t}(\bar{x},\bar{t}) \, |F(\bar{x},\bar{t}) - F(\bar{y},\bar{t})|^2 \\ 
&- H(\bar{x},\bar{t}) \, \Phi(\bar{x},\bar{t}) \, \langle F(\bar{x},\bar{t}) - F(\bar{y},\bar{t}),\nu(\bar{x},\bar{t}) \rangle \\ 
&- H(\bar{x},\bar{t}) + H(\bar{y},\bar{t}) \\ 
&+ \sum_{i=1}^n \frac{\partial H}{\partial x_i}(\bar{x},\bar{t}) \, \Big \langle F(\bar{x},\bar{t}) - F(\bar{y},\bar{t}),\frac{\partial F}{\partial x_i}(\bar{x},\bar{t}) \Big \rangle \\ 
&= \frac{1}{2} \, \bigg ( \frac{\partial \Phi}{\partial t}(\bar{x},\bar{t}) + H(\bar{x},\bar{t}) \, \Phi(\bar{x},\bar{t})^2 \\ 
&\hspace{20mm} - \sum_{i=1}^n \frac{1}{\Phi(\bar{x},\bar{t})+\lambda_i(\bar{x},\bar{t})} \, \frac{\partial \Phi}{\partial x_i}(\bar{x},\bar{t}) \, \frac{\partial H}{\partial x_i}(\bar{x},\bar{t}) \bigg ) \, |F(\bar{x},\bar{t}) - F(\bar{y},\bar{t})|^2 \\ 
&- H(\bar{x},\bar{t}) + H(\bar{y},\bar{t}). 
\end{align*} 
Note that the term $H(\bar{y},\bar{t})$ is nonnegative. In the next step, we multiply both sides by $\frac{2}{|F(\bar{x},\bar{t})-F(\bar{y},\bar{t})|^2}$. Using the identity 
\begin{align*} 
\frac{1}{|F(\bar{x},\bar{t})-F(\bar{y},\bar{t})|^2} &= \frac{\langle F(\bar{x},\bar{t})-F(\bar{y},\bar{t}),\nu(\bar{x},\bar{t}) \rangle^2}{|F(\bar{x},\bar{t})-F(\bar{y},\bar{t})|^4} + \sum_{i=1}^n \frac{\langle F(\bar{x},\bar{t})-F(\bar{y},\bar{t}),\frac{\partial F}{\partial x_i}(\bar{x},\bar{t}) \rangle^2}{|F(\bar{x},\bar{t})-F(\bar{y},\bar{t})|^4} \\ 
&= \frac{1}{4} \, \bigg ( \Phi(\bar{x},\bar{t})^2 + \sum_{i=1}^n \frac{1}{(\Phi(\bar{x},\bar{t})+\lambda_i(\bar{x},\bar{t}))^2} \, \Big ( \frac{\partial \Phi}{\partial x_i}(\bar{x},\bar{t}) \Big )^2 \bigg ), 
\end{align*} 
we deduce that 
\begin{align*} 
&\frac{2}{|F(\bar{x},\bar{t})-F(\bar{y},\bar{t})|^2} \, \frac{\partial W}{\partial t}(\bar{x},\bar{y},\bar{t}) \\ 
&\geq \frac{\partial \Phi}{\partial t}(\bar{x},\bar{t}) + H(\bar{x},\bar{t}) \, \Phi(\bar{x},\bar{t})^2 - \sum_{i=1}^n \frac{1}{\Phi(\bar{x},\bar{t})+\lambda_i(\bar{x},\bar{t})} \, \frac{\partial \Phi}{\partial x_i}(\bar{x},\bar{t}) \, \frac{\partial H}{\partial x_i}(\bar{x},\bar{t}) \\ 
&- \frac{1}{2} \, H(\bar{x},\bar{t}) \, \bigg ( \Phi(\bar{x},\bar{t})^2 + \sum_{i=1}^n \frac{1}{(\Phi(\bar{x},\bar{t})+\lambda_i(\bar{x},\bar{t}))^2} \, \Big ( \frac{\partial \Phi}{\partial x_i}(\bar{x},\bar{t}) \Big )^2 \bigg ) \\ 
&= \frac{\partial \Phi}{\partial t}(\bar{x},\bar{t}) + \frac{1}{2} \, H(\bar{x},\bar{t}) \, \Phi(\bar{x},\bar{t})^2 - \sum_{i=1}^n \frac{1}{\Phi(\bar{x},\bar{t})+\lambda_i(\bar{x},\bar{t})} \, \frac{\partial \Phi}{\partial x_i}(\bar{x},\bar{t}) \, \frac{\partial H}{\partial x_i}(\bar{x},\bar{t}) \\ 
&- \frac{1}{2} \, H(\bar{x},\bar{t}) \, \sum_{i=1}^n \frac{1}{(\Phi(\bar{x},\bar{t})+\lambda_i(\bar{x},\bar{t}))^2} \, \Big ( \frac{\partial \Phi}{\partial x_i}(\bar{x},\bar{t}) \Big )^2. 
\end{align*} 
Since $\frac{\partial W}{\partial t}(\bar{x},\bar{y},\bar{t}) \leq 0$, the assertion follows. \\

As above, we can show that 
\[\frac{\partial \rho}{\partial t} - \Delta \rho - |A|^2 \, \rho + \sum_{i=1}^n \frac{2}{\rho+\lambda_i} \, (D_i \rho)^2 \leq 0\] 
on the set $\{\rho > 0\} \cap \{\rho+\lambda_1 > 0\}$, where $\Delta \rho$ is interpreted in the sense of distributions. Moreover, Proposition \ref{aux.2} implies that 
\[\frac{\partial \rho}{\partial t} - \sum_{i=1}^n \frac{1}{\rho+\lambda_i} \, D_i \rho \, D_i H - \frac{1}{2} \, \sum_{i=1}^n \frac{H}{(\rho+\lambda_i)^2} \, (D_i \rho)^2 \leq 0\] 
almost everywhere on the set $\{\rho > 0\} \cap \{\rho+\lambda_1 > 0\}$. \\

Suppose now that a number $\delta>0$ is given. By the convexity estimate of Huisken and Sinestrari \cite{Huisken-Sinestrari2}, we can find a constant $K_0$, depending only on $\delta$ and the initial data, such that 
\[\lambda_1 \geq -\frac{\delta}{2} \, H - K_0 \, \min \{H,1\}.\] 
For each $\sigma \in (0,\frac{1}{2})$, we put 
\[g_\sigma = H^{\sigma-1} \, (\rho - \delta \, H) - K_0\] 
and 
\[g_{\sigma,+} = \max \{g_\sigma,0\}.\] 
On the set $\{g_\sigma \geq 0\}$, we have 
\[\rho \geq \delta \, H + K_0 \, H^{1-\sigma} \geq \delta \, H + K_0 \, \min \{H,1\},\] 
hence 
\[\rho + \lambda_1 \geq \frac{\delta}{2} \, H.\] 
In particular, we have $\{g_\sigma \geq 0\} \subset \{\rho > 0\} \cap \{\rho+\lambda_1 > 0\}$. \\

\begin{proposition}
\label{Lp.bound.2}
Given any $\delta>0$, we can find a positive constant $c_0$, depending only on $\delta$ and the initial data, with the following property: if $p \geq \frac{1}{c_0}$ and $\sigma \leq c_0 \, p^{-\frac{1}{2}}$, then we have 
\[\frac{d}{dt} \bigg ( \int_{M_t} g_{\sigma,+}^p \bigg ) \leq \sigma \, p \, K_0^p \int_{M_t} |A|^2\] 
for almost all $t$. Here, $C$ is a positive constant that depends only on $\delta$ and the initial data, but not on $\sigma$ and $p$.
\end{proposition}

\textbf{Proof.} 
For abbreviation, we define a function $\omega$ by 
\begin{align*} 
\omega 
&= \Delta \rho - \sum_{i=1}^n \frac{2}{\rho+\lambda_i} \, (D_i \rho)^2 \\ 
&- \sum_i \frac{1}{\rho+\lambda_i} \, D_i \rho \, D_i H - \frac{1}{2} \, \sum_{i=1}^n \frac{H}{(\rho+\lambda_i)^2} \, (D_i \rho)^2. 
\end{align*}
The function $\rho$ satisfies 
\[\frac{\partial \rho}{\partial t} - \Delta \rho - |A|^2 \, \rho + \sum_{i=1}^n \frac{2}{\rho+\lambda_i} \, (D_i \rho)^2 \leq -\max \{\omega + |A|^2 \, \rho,0\}\] 
on the set $\{\rho > 0\} \cap \{\rho+\lambda_1 > 0\}$. From this, we deduce that 
\begin{align*} 
&\frac{\partial}{\partial t} g_\sigma - \Delta g_\sigma - 2 \, (1-\sigma) \, \Big \langle \frac{\nabla H}{H},\nabla g_\sigma \Big \rangle \\ 
&+ 2 \, \sum_{i=1}^n \frac{H^{\sigma-1}}{\rho+\lambda_i} \, (D_i \rho)^2 - \sigma \, |A|^2 \, (g_\sigma + K_0) \\ 
&\leq -H^{\sigma-1} \, \max \{\omega + |A|^2 \, \rho,0\} - \sigma \, (1-\sigma) \, H^{\sigma-3} \, (\rho - \delta \, H) \, |\nabla H|^2 
\end{align*} 
on the set $\{g_\sigma \geq 0\}$. Note that $g_\sigma \leq H^{\sigma-1} \, \rho$ by definition of $g_\sigma$. Since $\sigma \in (0,\frac{1}{2})$, we have $2\sigma \, \frac{g_\sigma}{\rho} \leq H^{\sigma-1}$ at each point on the surface. This implies 
\begin{align*} 
&\frac{\partial}{\partial t} g_\sigma - \Delta g_\sigma - 2 \, (1-\sigma) \, \Big \langle \frac{\nabla H}{H},\nabla g_\sigma \Big \rangle + 2 \, \sum_{i=1}^n \frac{H^{\sigma-1}}{\rho+\lambda_i} \, (D_i \rho)^2 \\ 
&\leq -H^{\sigma-1} \, \max \{\omega + |A|^2 \, \rho,0\} + \sigma \, |A|^2 \, (g_\sigma + K_0) \\ 
&\leq -2\sigma \, \frac{g_\sigma}{\rho} \, (\omega + |A|^2 \, \rho) + \sigma \, |A|^2 \, (g_\sigma + K_0) \\ 
&= -2\sigma \, \frac{g_\sigma}{\rho} \, \omega + \sigma \, |A|^2 \, (K_0-g_\sigma) 
\end{align*} 
on the set $\{g_\sigma \geq 0\}$. Therefore, we have 
\begin{align*} 
&\frac{d}{dt} \bigg ( \int_{M_t} g_{\sigma,+}^p \bigg ) \\ 
&\leq -p(p-1) \int_{M_t} g_{\sigma,+}^{p-2} \, |\nabla g_\sigma|^2 + 2 \, (1-\sigma) \, p \int_{M_t} g_{\sigma,+}^{p-1} \, \Big \langle \frac{\nabla H}{H},\nabla g_\sigma \Big \rangle \\ 
&- 2p \int_{M_t} \sum_{i=1}^n \frac{g_{\sigma,+}^{p-1} \, H^{\sigma-1}}{\rho+\lambda_i} \, (D_i \rho)^2 - 2\sigma \, p \int_{M_t} \frac{g_{\sigma,+}^p}{\rho} \, \omega \\ 
&+ \sigma \, p \int_{M_t} g_{\sigma,+}^{p-1} \, (K_0-g_\sigma) \, |A|^2. 
\end{align*} 
Integration by parts gives 
\begin{align*} 
-\int_{M_t} \frac{g_{\sigma,+}^p}{\rho} \, \omega &\leq C \, p \int_{M_t} \frac{g_{\sigma,+}^{p-1}}{H} \, |\nabla \rho| \, |\nabla g_\sigma| \\ 
&+ C \int_{M_t} \frac{g_{\sigma,+}^p}{H^2} \, |\nabla \rho| \, |\nabla H| + C \int_{M_t} \frac{g_{\sigma,+}^p}{H^2} \, |\nabla \rho|^2, 
\end{align*} 
where $C$ is a positive constant that depends only on $\delta$ and the initial data. Putting these facts together, we obtain 
\begin{align*} 
&\frac{d}{dt} \bigg ( \int_{M_t} g_{\sigma,+}^p \bigg ) \\ 
&\leq -p(p-1) \int_{M_t} g_{\sigma,+}^{p-2} \, |\nabla g_\sigma|^2 + 2 \, (1-\sigma) \, p \int_{M_t} g_{\sigma,+}^{p-1} \, \Big \langle \frac{\nabla H}{H},\nabla g_\sigma \Big \rangle \\ 
&- 2p \int_{M_t} \sum_{i=1}^n \frac{g_{\sigma,+}^{p-1} \, H^{\sigma-1}}{\rho+\lambda_i} \, (D_i \rho)^2 + C \, \sigma \, p^2 \int_{M_t} \frac{g_{\sigma,+}^{p-1}}{H} \, |\nabla \rho| \, |\nabla g_\sigma| \\ 
&+ C \, \sigma \, p \int_{M_t} \frac{g_{\sigma,+}^p}{H^2} \, |\nabla \rho| \, |\nabla H| + C \, \sigma \, p \int_{M_t} \frac{g_{\sigma,+}^p}{H^2} \, |\nabla \rho|^2 \\ 
&+ \sigma \, p \, K_0^p \int_{M_t} |A|^2, 
\end{align*} 
where $C$ is a positive constant that depends only on $\delta$ and the initial data. Using the identity 
\[\frac{\nabla H}{H} = \frac{\nabla \rho - H^{1-\sigma} \, \nabla g_\sigma}{(1-\sigma) \, \rho + \sigma \, \delta \, H},\] 
we obtain 
\[\Big \langle \frac{\nabla H}{H},\nabla g_\sigma \Big \rangle \leq \frac{\langle \nabla \rho,\nabla g_\sigma \rangle}{(1-\sigma) \, \rho + \sigma \, \delta \, H} \leq C \, H^{-1} \, |\nabla \rho| \, |\nabla g_\sigma|\] 
and 
\[\frac{|\nabla H|}{H} \leq C \, \frac{|\nabla \rho|}{H} + C \, \frac{|\nabla g_\sigma|}{g_{\sigma,+}}.\] 
This gives 
\begin{align*} 
&\frac{d}{dt} \bigg ( \int_{M_t} g_{\sigma,+}^p \bigg ) \\ 
&\leq -p(p-1) \int_{M_t} g_{\sigma,+}^{p-2} \, |\nabla g_\sigma|^2 - 2p \int_{M_t} \sum_{i=1}^n \frac{g_{\sigma,+}^{p-1} \, H^{\sigma-1}}{\rho+\lambda_i} \, (D_i \rho)^2 \\ 
&+ C \, (p + \sigma \, p^2) \int_{M_t} \frac{g_{\sigma,+}^{p-1}}{H} \, |\nabla \rho| \, |\nabla g_\sigma| + C \, \sigma \, p \int_{M_t} \frac{g_{\sigma,+}^p}{H^2} \, |\nabla \rho|^2 \\ 
&+ \sigma \, p \, K_0^p \int_{M_t} |A|^2, 
\end{align*} 
where $C$ is a positive constant that depends only on $\delta$ and the initial data.

Consequently, there exists a positive constant $c_0$, depending only on $\delta$ and the initial data, with the following property: if $p \geq \frac{1}{c_0}$ and $\sigma \leq c_0 \, p^{-\frac{1}{2}}$, then we have 
\[\frac{d}{dt} \bigg ( \int_{M_t} g_{\sigma,+}^p \bigg ) \leq \sigma \, p \, K_0^p \int_{M_t} |A|^2.\] 
This completes the proof of Proposition \ref{Lp.bound.2}. \\

Moreover, we have the following estimate:

\begin{proposition}
Let $g_{\sigma,k} = H^{\sigma-1} \, (\rho - \delta \, H) - k$ and $g_{\sigma,k,+} = \max \{g_{\sigma,k},0\}$. Then 
\begin{align*} 
\frac{d}{dt} \bigg ( \int_{M_t} g_{\sigma,k,+}^p \bigg ) 
&\leq -\frac{1}{2} \, p(p-1) \int_{M_t} g_{\sigma,k,+}^{p-2} \, |\nabla g_{\sigma,k}|^2 \\ 
&+ \sigma \, p \int_{M_t} |A|^2 \, g_{\sigma,k,+}^{p-1} \, (g_{\sigma,k}+k) 
\end{align*}
if $k \geq K_0$ and $p$ is sufficiently large.
\end{proposition} 

\textbf{Proof.} 
Assume that $k \geq K_0$. Using the inequality 
\begin{align*} 
\frac{\partial}{\partial t} g_{\sigma,k} 
&\leq \Delta g_{\sigma,k} + 2 \, (1-\sigma) \, \Big \langle \frac{\nabla H}{H},\nabla g_{\sigma,k} \Big \rangle \\ 
&- 2 \, \sum_{i=1}^n \frac{H^{\sigma-1}}{\rho+\lambda_i} \, (D_i \rho)^2 + \sigma \, |A|^2 \, (g_{\sigma,k} + k), 
\end{align*} 
we obtain 
\begin{align*}
&\frac{d}{dt} \bigg ( \int_{M_t} g_{\sigma,k,+}^p \bigg ) \\ 
&\leq -p(p-1) \int_{M_t} g_{\sigma,k,+}^{p-2} \, |\nabla g_{\sigma,k}|^2 + 2 \, (1-\sigma) \, p \int_{M_t} g_{\sigma,k,+}^{p-1} \, \Big \langle \frac{\nabla H}{H},\nabla g_{\sigma,k} \Big \rangle \\ 
&- 2p \int_{M_t} \sum_{i=1}^n \frac{g_{\sigma,k,+}^{p-1} \, H^{\sigma-1}}{\rho+\lambda_i} \, (D_i \rho)^2 + \sigma \, p \int_{M_t} |A|^2 \, g_{\sigma,k,+}^{p-1} \, (g_{\sigma,k}+k). 
\end{align*} 
As above, we have 
\[\Big \langle \frac{\nabla H}{H},\nabla g_{\sigma,k} \Big \rangle \leq \frac{\langle \nabla \rho,\nabla g_{\sigma,k} \rangle}{(1-\sigma) \, \rho + \sigma \, \delta \, H} \leq C \, H^{-1} \, |\nabla \rho| \, |\nabla g_{\sigma,k}|,\] 
where $C$ depends only on $\delta$ and the initial data. This gives 
\begin{align*}
&\frac{d}{dt} \bigg ( \int_{M_t} g_{\sigma,k,+}^p \bigg ) \\ 
&\leq -p(p-1) \int_{M_t} g_{\sigma,k,+}^{p-2} \, |\nabla g_{\sigma,k}|^2 + C \, p \int_{M_t} \frac{g_{\sigma,k,+}^{p-1}}{H} \, |\nabla \rho| \, |\nabla f_{\sigma,k}| \\ 
&- 2p \int_{M_t} \sum_{i=1}^n \frac{g_{\sigma,k,+}^{p-1} \, H^{\sigma-1}}{\rho+\lambda_i} \, (D_i \rho)^2 + \sigma \, p \int_{M_t} |A|^2 \, g_{\sigma,k,+}^{p-1} \, (g_{\sigma,k}+k). 
\end{align*} From this, the assertion follows. \\

As above, we can use Stampacchia iteration to show that $g_\sigma \leq k$, where $\sigma$ and $k$ depend only on $\delta$ and the initial data. From this, we deduce that 
\[\rho \leq 2\delta \, H + B,\] 
where $B$ depends only on $\delta$ and the initial data. This completes the proof of Theorem \ref{main.theorem.2}.

\end{document}